\documentclass[12pt,a4paper]{amsart}

\usepackage{amsthm}
\usepackage{amsmath}
\usepackage{amssymb}
\usepackage{amsfonts}
\usepackage{amscd}
\usepackage{palatino}
\usepackage[mathscr]{eucal}
\usepackage{epsfig}
\usepackage{verbatim}
\usepackage{latexsym}
\usepackage{graphics}
\usepackage{color}
\usepackage{a4}
\usepackage{graphicx}
\usepackage{appendix}

\begin{document}

\newcommand{\fr}{{\rm Frac}}

\title{The number of rational numbers determined by
large sets of integers}

\author{Javier Cilleruelo}
\address{Instituto de Ciencias Matem\'{a}ticas (CSIC-UAM-UC3M-UCM) and Departamento de
Matem\'{a}ticas,
Universidad Aut\'{o}noma de Madrid,
Madrid-28049, Spain}
\email{franciscojavier.cilleruelo@uam.es}

\author{D.S. Ramana}
\address{ Harish-Chandra Research Institute, Jhunsi, Allahabad -211 019, India.}
\email{suri@hri.res.in }

\author{Olivier Ramar\'{e}}
\address{Laboratoire Paul Painlev\'{e},
Universit\'{e} Lille 1,
59655 Villeneuve d'Ascq Cedex, France}
\email{ramare@math.univ-lille1.fr}

\makeatletter{\renewcommand*{\@makefnmark}{}
\makeatother}

\begin{abstract} When  $A$ and $B$ are subsets of the integers in $[1,X]$ and $[1,Y]$ respectively,
with $|A| \geq \alpha X$ and $|B| \geq \beta X$, we show that the
number of rational numbers expressible as $a/b$ with $(a,b)$ in $A
\times B$ is $\gg (\alpha \beta)^{1+\epsilon}XY$ for any $\epsilon >
0$, where the implied constant depends on $\epsilon$ alone. We then
construct examples that show that this bound cannot in general be
improved to $\gg \alpha \beta XY$. We also resolve the natural 
generalisation of our problem to arbitrary subsets $C$ of the
integer points in $[1,X] \times [1,Y]$. Finally,
we apply our results to answer a question of S\'ark\"ozy concerning
the differences of consecutive terms of the product sequence of a
given integer sequence.\end{abstract}

\maketitle

\newcounter{thn}
\renewenvironment{th}
{\addtocounter{thn}{1}
 \em
     \noindent
     {\sc Theorem \arabic{scn}.\arabic{thn}. ---}}{}

\newcounter{len}
\renewenvironment{le}
{\addtocounter{len}{1}
 \em
     \vspace{2mm}
     \noindent
     {\sc Lemma \arabic{scn}.\arabic{len}. ---}}{}

\newcounter{prn}
\newenvironment{prop}
{\addtocounter{prn}{1}
  \em
      \vspace{2mm}
      \noindent
      {\sc Proposition \arabic{scn}.\arabic{prn}. ---}}{}

\newcounter{corn}
\newenvironment{cor}
  {\em
     \addtocounter{corn}{1}
     \vspace{2mm}
     \noindent
        {\sc Corollary \arabic{scn}.\arabic{corn}. ---}}{\vspace{1mm}}

\newcounter{scn}
\newenvironment{nsc}[1]
              {\addtocounter{scn}{1}
              \vspace{3mm}
                     \begin{center}
                     {\sc { \arabic{scn}. {#1}}}
                     \end{center}
              \setcounter{equation}{0}
              \setcounter{prn}{0}
              \setcounter{len}{0}
              \setcounter{thn}{0}
              \setcounter{corn}{0}
              \vspace{3mm}}{}

\begin{nsc}{Introduction}
When $A$ and $B$ are intervals in the integers in $[1,X]$ and
$[1,Y]$ respectively, satisfying $|A| \geq \alpha X$ and $|B| \geq
\beta Y$, where $X$, $Y$ real numbers $\geq 1$, $\alpha$,  $\beta$ are real numbers in $(0,1]$, a standard application of the M\"{o}bius
inversion formula shows that the number of rational numbers $a/b$
with $(a,b)$ in $A \times B$ is $\gg \alpha \beta XY$. 

\vspace{2mm} \noindent Our purpose is to investigate what might be
deduced when in place of {\em intervals} we consider {\em arbitrary}
subsets $A$ and $B$ of the integers in $[1,X]$ and $[1,Y]$
respectively with $|A| \geq \alpha X$ and $|B| \geq \beta Y$. When
$A$ and $B$ are not intervals, it may happen that an abnormally
large number of elements of these sets are multiples of certain
integers, determining which in general is not easy. 
Nevertheless, since the sets under consideration are
large, popular heuristics suggest that a
non-trivial conclusion should still be accessible. What is pleasing
is that we in fact have the following theorem, which is our principal conclusion. In the
statement of this theorem and thereafter we write  $A/B$ to denote
the subset of ${\bf Q}$ consisting of all rational numbers
expressible as $a/b$ with $(a,b)$ in $A \times B$ for any $A$ and
$B$ subsets of the integers $\geq 1$.

\vspace{3mm}
\begin{th}\label{main}
Let $\alpha$ and $\beta$ be real numbers in $(0,1]$ and $X$ and $Y$
real numbers $\geq 1$. When $A$ and $B$ are subsets of the integers
in $[1,X]$ and $[1,Y]$ respectively, with $|A| \geq
\alpha X$ and $|B| \geq \beta Y$ we have $|A/B| \gg
(\alpha\beta)^{1+\epsilon}XY$ for any $\epsilon > 0$, where the
implied constant depends on $\epsilon$ alone.
\end{th}

\vspace{2mm} \noindent Deferring the detailed proof of Theorem 1.1 to Section~ 2, let us summarize our argument with the aid of the following notation. For any integer  $d \geq 1$,  $A$ and $B$ subsets of the integers $\geq 1$, we write ${\mathcal
M}(A,B,d)$ to denote the subset of $A \times B$ consisting of all $(a,b)$ in $A \times B$ with $\gcd(a,b)=d$.
We show in Proposition 2.1 that for $A$ and $B$ as in Theorem 1.1 we have $\sup_{d \geq 1} |{\mathcal M}(A,B,d)| \geq \frac{1}{8} (\alpha \beta)^{2} XY$. Starting from this  initial bound we then obtain $\sup_{d \geq 1} |{\mathcal M}(A,B,d)| \gg (\alpha \beta)^{1 +\epsilon} XY$ by a bootstrapping argument. Theorem 1.1 follows immediately from this last bound, since for any integer $d \geq 1$ we have $a/b \neq a_1/b_1$ for any two points $(a,b)$ and $(a_1,b_1)$ of
${\mathcal M}(A,B,d)$, and therefore $|A/B|
\geq \sup_{d \geq 1} |{\mathcal M}(A,B,d)|$.

\vspace{2mm} \noindent
We supplement Theorem 1.1 with the following result which shows that the bound provided by Theorem 1.1 cannot be replaced with $|A/B| \gg \alpha\beta XY$. This bound, as we have already remarked, holds when $A$ and $B$ are intervals.

\vspace{3mm}
\begin{th}
For any $\epsilon >0$, there exists $\alpha>0$ such that for all
sufficiently large $X$ there exists a subset $A$ of the integers in
$[1,X]$ satisfying $|A|\ge \alpha X$ and $|A/A|<\epsilon
\alpha^2X^2$.
\end{th}

\vspace{2mm}
\noindent
We prove Theorem 1.2 in Section 3.
Our method depends on the observation that
for any $\epsilon > 0$ and any set of prime numbers ${\mathcal P}$ with  $|{\mathcal P}|$ sufficiently large, we have $|
S({\mathcal P})/ S({\mathcal P})| \leq
\epsilon |S({\mathcal P})|^2$
, where $S({\mathcal P})$ is the set of squarefree
integers $d$ formed from the primes in the subsets of ${\mathcal P}$ containing about half the primes in ${\mathcal P}$.
 By means of this observation we deduce that, for a suitable ${\mathcal P}$, the set of multiples of the elements of $S({\mathcal P})$ in $[1,X]$, meets the conditions of Theorem 1.2.

\vspace{2mm}
\noindent
The questions answered by the above theorems may be viewed as particular cases of a more general problem namely, for $X$ and $Y$ real numbers $\geq 1$ and $\gamma$ in $(0,1]$, given a subset $C$ of the integer points in $[1,X] \times [1,Y]$ satisfying $|C| \geq \gamma XY$, to determine in terms of $\gamma$, $X$ and $Y$ an optimal lower bound for ${\rm Frac}(C)$, the number of rational numbers $a/b$ with $(a,b)$ in $C$. Plainly, the above theorems take up the special case when $C$ is of the form $A \times B$, that is, when $C$ is equal to the product of its projections onto the co-ordinate axes.

\vspace{2mm}
\noindent
It turns out, however, that aforementioned general problem is somewhat easily resolved. In effect, the method of Proposition 2.1 generalizes without additional effort to give the bound $|{\rm Frac}(C)| \geq \frac{1}{8}\gamma^2 XY$ and, interestingly, this bound is in fact optimal upto the constant $\frac{1}{8}$. More precisely, we have the following theorem.

\vspace{3mm}
\begin{th}
For any $\gamma$ in $(0,1]$ and all sufficiently large $X$ and $Y$ there exists a subset $C$ of the integer points in $[1,X] \times [1,Y]$ satisfying $|C| \geq \frac{\gamma}{8} XY$ and $|\fr(C)| \leq \frac{\gamma^2}{2}XY$.
\end{th}

\vspace{2mm}
\noindent
We prove Theorem 1.3 at the end of Section 3 by explicitly describing sets $C$ that satisfy the conditions of this theorem. Such sets are in general far from being of the form $A \times B$, which is only natural on account of Theorem 1.1. Indeed, our bootstrapping argument for Theorem 1.1 depends crucially on the fact that this theorem is, from the more general point view, about sets $C$ which are of the form $A \times B$.

\vspace{2mm} \noindent We conclude this note with Section 4 where we
 apply Theorem 1.1 to obtain a near-optimal answer to
the following question of S\'ark\"ozy. When ${\mathcal A}$, ${\mathcal B}$ are
sequences of integers, let ${\mathcal A}.{\mathcal B}$ be the sequence whose
terms are the integers of the form $ab$, for some  $a\in
{\mathcal A}$, $b\in {\mathcal B}$. Then S\'ark\"ozy \cite{sar} asks if it is true that for
any $\alpha
> 0$ and ${\mathcal A}$ such that the lower asymptotic density
$\underline{d}({\mathcal A}) > \alpha$ there is a $c(\alpha)$ such that there
are infinitely many pairs of consecutive terms of ${\mathcal A}.{\mathcal A}$ the difference
between which is bounded by $c(\alpha)$.

\vspace{2mm} \noindent Berczi \cite{ber} responded to the
aforementioned question of S\'ark\"ozy by showing that the minimum of
the differences between consecutive terms of ${\mathcal A}.{\mathcal A}$ is $\ll
\frac{1}{\alpha^4}$, where $\alpha = \underline{d}({\mathcal A})$. Sandor
\cite{san} subsequently improved this by showing that this minimum
is in fact $\ll \frac{1}{\alpha^3}$, with $\alpha$ now the upper
asymptotic density $\overline{d}({\mathcal A})$ of ${\mathcal A}$. Cilleruelo and Le
\cite{ch} obtained the same bound when $\alpha$ is the upper Banach
density of ${\mathcal A}$ and showed that this is the best possible bound for
this density. The following result improves upon and
generalizes Sandor's conclusion.

\vspace{3mm}
\begin{th}
Let $\alpha$ and $\beta$ be real numbers in $(0,1]$ and let
$\epsilon$ be $> 0$.  When ${\mathcal A}$ and ${\mathcal B}$ are infinite sequences of
integers with upper asymptotic densities $\alpha$ and $\beta$
respectively, there are infinitely many pairs of
consecutive terms of the product sequence ${\mathcal A}.{\mathcal B}$ the difference
between which is $\ll \frac{1}{(\alpha \beta)^{1+\epsilon}}$, where
the implied constant depends  on $\epsilon$ alone.
\end{th}

\vspace{2mm} \noindent When ${\mathcal A}$ and ${\mathcal B}$ are the sequences of
multiples of the integers $h$ and $k$ respectively, the difference
between any two consecutive terms of the sequence ${\mathcal A}.{\mathcal B}$ is $\geq
hk$. Since we have  $\overline{d}({\mathcal A})= \frac{1}{h}$ and
$\overline{d}({\mathcal B}) = \frac{1}{k}$, we see that the conclusion of
Theorem 1.4 is optimal up to a factor $\frac
1{(\alpha \beta)^{\epsilon}}$.

\vspace{2mm}
\noindent
Throughout this note, $X$, $Y$ shall denote real numbers $\geq 1$ and $\alpha$, $\beta$, $\gamma$ 
real numbers in $(0,1]$. Also, the letter $p$ shall denote a prime number. When $I$ and $J$ are subsets of a given set, $I \setminus J$ shall denote the set of elements of $I$ that are not in $J$. In addition to  
the notation introduced so far, we shall write $A_d$ to denote the subset of a set of
integers $A$ consisting of all multiples of $d$ in $A$ for any integer $d$. Finally, if $B = \{b\}$ with
$b \geq 1$, we simply write $A/b$ in place of $A/B$, by an abuse of notation.

\end{nsc}

\begin{nsc}{Proof of the Bound}

\noindent
Let $A$ and $B$ be finite subsets of the integers $\geq 1$.
Then the family of subsets ${\mathcal M}(A,B,d)$ of $A \times B$, with $d$
varying over the integers $\geq 1$, is a partition of $A \times B$. Consequently, we have

\vspace{-2mm}
\begin{equation}
\label{f1}
|A \times B| \; = \; \sum_{d \geq 1} |{\mathcal M}(A,B,d)| \; .
\end{equation}

\noindent When $A$ and $B$ are contained in $[1,X]$ and $[1,Y]$
respectively, we have $|A_d| \leq X/d$ and $|B_d| \leq
Y/d$, for any $d \geq 1$. Since ${\mathcal M}(A,B,d)$ is contained
in $A_d \times B_d$, we then obtain $|{\mathcal M}(A,B,d)| \leq
|A_d||B_d| \leq \frac{XY}{d^2}$, for all $d \geq 1$.

\begin{prop}
When $A$ and $B$ are subsets of the integers in the intervals $[1,X]$ and $[1,Y]$
respectively, with $|A| \geq \alpha X$ and $|B| \geq \beta Y$, we have $\sup_{d \geq 1} |{\mathcal M}(A,B,d)| \geq \frac{(\alpha \beta)^2 XY}{8}$.
\end{prop}

\vspace{3mm}
\noindent
{\sc proof.---} We adapt an argument from \cite{ch}. From (\ref{f1}) we have for any integer $T \geq 1$ that

\vspace{-2mm}
\begin{equation}
\label{f2}
|A \times B| \; = \;
\sum_{1 \leq d \leq T} |{\mathcal M}(A,B,d)| + \sum_{T < d} |{\mathcal M}(A,B,d)| \; \leq \;
\sum_{1 \leq d \leq T} |{\mathcal M}(A,B,d)| + \frac{XY}{T} \; ,
\end{equation}

\noindent
where the last inequality follows from $\sum_{T < d} |{\mathcal M}(A,B,d)| \leq \sum_{T < d} \frac{XY}{d^2} \leq \frac{XY}{T}$.
Since $|A \times B| \geq \alpha \beta XY$ we conclude
from (\ref{f2}) that

\vspace{-2mm}
\begin{equation}
\label{f3}
\sup_{d \geq 1} |{\mathcal M}(A,B,d)|
\; \geq \;
\frac{1}{T} \sum_{1 \leq d \leq T} |{\mathcal M}(A,B,d)|
\; \geq \;
\left(\frac{\alpha\beta - \frac{1}{T}}{T}\right)XY \;
\end{equation}

\noindent
for any integer $T \geq 1$. Since $2 > \alpha \beta$, the interval $[\frac{2}{\alpha\beta}, \frac{4}{\alpha\beta}]$
contains an integer $\geq 1$. The proposition now follows on setting $T$ in (\ref{f3}) to be any such integer.

\vspace{3mm} \noindent {\sc Definition 2.1---} We call a real number
$\delta$ an {\em admissible exponent} if there exists a real number
$C >0$ such that for any $\alpha$, $\beta$ real numbers in $(0,1]$,
any $X$, $Y$ real numbers $\geq 1$ and any subsets $A$ and $B$ of
the integers in $[1,X]$ and $[1,Y]$ with $|A| \geq \alpha X$ and $|B| \geq \beta Y$, we have $\sup_{d \geq 1} |{\mathcal M}(A,B,d)| \geq
C(\alpha\beta)^{\delta} XY$. We call a $C$ satisfying
these conditions a {\em constant associated to} the admissible
exponent $\delta$.

\vspace{2mm} \noindent
Proposition 2.1 says that $\delta=2$ is an
admissible exponent. Proposition 2.2 will allow us to
conclude that every $\delta > 1$ is an admissible exponent.
The following lemma prepares us for an application of H\"{o}lder's inequality within the proof of Proposition 2.2.

\vspace{2mm}
\noindent
For any integer $n\geq 1$ let $\tau(n)$ denote, as usual, the number of integers $\geq 1$ that divide $n$.
 When $D$ is an integer $\geq 1$ we write $\tau_D(n)$ to denote the number of divisors $d$ of $n$ satisfying the condition $p|d \implies p \leq D$ for any prime number $p$.

\vspace{2mm}
\begin{le}
When $q$ is an integer $\geq 0$ there is a real number $c(q) >0$ such that for all real numbers $X \geq 1$ and integers $D\geq 1$ we have

\vspace{-2mm}
\begin{equation}
\label{f4} \sum_{1 \leq n \leq X} \tau_{D}(n)^q \; \leq \;
c(q) DX \; ,
\end{equation}
\end{le}

\vspace{2mm}
\noindent
{\sc Proof.---} In effect, we have

\vspace{-2mm}
\begin{equation}
\label{6}
\sum_{1 \leq n \leq X}\tau_{D}(n)^q \ll X(\log 2D)^{2^q}
 \, \ll \, (2^q!)\,D X \; ,
\end{equation}

\noindent
where the implied constants are absolute. Plainly, the second inequality results from the elementary inequality $(\log t)^n \leq n!\, t$ for $t \geq 1$. We now
prove the first inequality in (\ref{6}). Let us write ${\mathcal D}$ for the set of integers
$m$ satisfying the condition $p|m \implies p \leq D$. For any integer $n \geq 1$, let $k(n)$
be the largest of the divisors of $n$ lying in ${\mathcal D}$. We then have that

\begin{equation}
\label{7}
\sum_{1 \leq n \leq X}\tau_{D}(n)^q = \sum_{\stackrel{1 \leq m \leq X,}{m \in {\mathcal D}}} \tau(m)^q \sum_{\stackrel{1 \leq n \leq X,}{k(n) = m}} 1
\leq X \sum_{m \in {\mathcal D}} \frac{\tau(m)^q}{m} \;
\; ,
\end{equation}

\noindent
where we have used the upper bound $X/m$ for the number of integers $n$ in $[1,X]$ with $k(n)=m$. Let us write $S(q)$ for any integer $q \geq 0$ to denote the last sum in (\ref{7}). Since Merten's formula gives
$\prod_{1 \leq p \leq D} (1-\frac{1}{p}) \sim \frac{e^{-\gamma}}{\log D}$, with $\gamma$ here being Euler's constant, we have

\vspace{-2mm}
\begin{equation}
\label{9}
S(0) = \sum_{ m \in D} \frac{1}{m} \; = \;
\prod_{1 \leq p \leq D}\left(1 + \frac{1}{p} + \frac{1}{p^2} +\ldots\right)\; = \;
\prod_{1 \leq p \leq D} \left(1 - \frac{1}{p}\right)^{-1} \; \ll \; \log2D \; ,
\end{equation}

\noindent
where the implied constant is absolute. On noting that every divisor of an integer in ${\mathcal D}$ is again in ${\mathcal D}$ and using $\tau(dk) \leq \tau(d)\tau(k)$, valid for any integers $d$ and $k \geq 1$, we  obtain

\vspace{-2mm}
\begin{equation}
\label{8}
\sum_{m \in D} \frac{\tau(m)^q}{m} \;
\;= \;
\sum_{m \in D} \frac{\tau(m)^{q-1}}{m} \sum_{d | m} 1 \;
= \;
\sum_{(d,k) \in D\times D} \frac{\tau(dk)^{q-1}}{dk}
\; \leq \;
\left(\sum_{d \in D} \frac{\tau (d)^{q-1}}{d}\right)^2 \; .
\end{equation}

\noindent
In other words, $S(q) \leq S(q-1)^2$, for any $q \geq 1$. An induction on $q$ then shows that for any integer $q \geq 0$ we have $S(q) \leq S(0)^{2^q} \ll (\log D)^{2^q}$, where the implied constant is absolute. On combining this bound with (\ref{7}) we obtain the first inequality in (\ref{6}).

\vspace{3mm}
\begin{prop}
If $\delta>1$ is an admissible exponent then so is
$\frac{3\delta(1+1/q)-2}{2\delta-1}$ for every integer $q \geq 1$.
\end{prop}

\vspace{2mm} \noindent {\sc Proof.---} Let $q$ be a given integer
$\geq 1$ and, for the sake of conciseness, let us write
$\delta^{\prime}$ to denote $\frac{3\delta(1+1/q)-2}{2\delta-1}$,
 which is $> 1$ since $\delta > 1$.

\vspace{2mm}
\noindent
When $C$ is a constant associated to $\delta$, let us set $C^{\prime}$ to be the unique real number $> 0$ satisfying

\vspace{-2mm}
\begin{equation}
\label{f7.5} \frac{1}{8C^{\prime}} \; =\;
\left(\frac{C^{\prime}}{C}\right)^{\frac{1}{2(\delta-1)}}
8^{\frac{\delta}{\delta-1}}
(4c(q))^{\frac{\delta}{q(\delta-1)}} \; ,
\end{equation}

\vspace{2mm} \noindent where $c(q)$ is the implied constant in
(\ref{f4}) of Lemma 2.1. It is easily seen from (\ref{f7.5}) that by replacing
$C$ with a smaller constant associated to $\delta$ if necessary we
may assume that $\frac{1}{4} \geq C^{\prime}$.

\vspace{2mm}
\noindent
We shall show that $\delta^{\prime}$ is an admissible exponent with $C^{\prime}$ a constant associated to $\delta^{\prime}$.
Thus let $\alpha$, $\beta$ be real numbers in $(0,1]$ and $X$, $Y$ real numbers $\geq 1$. Also, let $A$ and $B$
 be any subsets of the integers in $[1,X]$ and $[1,Y]$ satisfying $|A| \geq \alpha X$ and $|B| \geq  \beta Y$. We shall show that

\vspace{-2mm}
\begin{equation}
\label{f8}
\sup_{d \geq 1} |{\mathcal M}(A,B,d)| \geq C^{\prime} {(\alpha \beta)}^{\delta^{\prime}}XY \; .
\end{equation}

\vspace{2mm}
\noindent
Replacing $\alpha$ and $\beta$ with $\alpha^{\prime}\geq \alpha$ and $\beta^{\prime}\geq \beta$ such that $ \alpha^{\prime}\leq |A| \leq 2\alpha^{\prime}$ and $\beta^{\prime}\leq |B| \leq 2\beta^{\prime}$ if necessary, we reduce to the case when $|A| \leq 2\alpha X$ and $|B| \leq 2\beta Y$.

\vspace{2mm}
\noindent
Let us first dispose of the possibility that an abnormally large number of the integers in $A$ and $B$ are multiples of a given integer.
Thus let $\alpha_d = |A_d|/X$ and $\beta_d = |B_d|/Y$, for any integer $d \geq 1$. Suppose that there exists an integer $d \geq 1$ such that

\vspace{-2mm}
\begin{equation}
\label{f9}
\alpha_d \beta_d \geq \left(\frac{C^{\prime}}{C}\right)^{\frac{1}{\delta}} (\alpha \beta)^{\frac{\delta^{\prime}}{\delta}} d^{\frac{2}{\delta}-2}\; .
\end{equation}

\noindent
Then $A_d$ and $B_d$ are both non-empty and therefore $X$ and $Y$ are both $\geq d$. Further,
the sets $A_d/d$ and $B_d/d$ are subsets of the integers in $[1,X/d]$ and $[1,Y/d]$.
Since $\delta$ is an admissible exponent, $C$ a constant associated to $\delta$, and
we have $|A_d/d|= (d\alpha_d)|X/d|$, $|B_d/d|= (d\beta_d)|X/d|$, there exists an integer $d^{\prime} \geq 1$ such that

\vspace{-2mm}
\begin{equation}
\label{f10}
|{\mathcal M}(A_d/d,B_d/d,d^{\prime})|\; \geq \; C(d^2\alpha_d\beta_d)^{\delta}\frac{XY}{d^2}\; \geq \; C^{\prime}(\alpha \beta)^{\delta^{\prime}} XY \;  ,
\end{equation}

\noindent
where the last inequality follows from (\ref{f9}). Since
$|{\mathcal M}(A_d/d,B_d/d,d^{\prime})|$ does not exceed $|{\mathcal M}(A,B,dd^{\prime})|$,
we obtain (\ref{f8}) from (\ref{f10}). We may therefore verify (\ref{f8}) assuming that for every integer $d \geq 1$ we have

\vspace{-2mm}
\begin{equation}
\label{f11} \alpha_d \beta_d <
\left(\frac{C^{\prime}}{C}\right)^{\frac{1}{ \delta}} (\alpha
\beta)^{\frac{\delta^{\prime}}{\delta}} d^{\frac{ 2}{\delta}-2}\; .
\end{equation}

\noindent
With the aid of (\ref{f11}) we shall in fact obtain a more precise conclusion than (\ref{f8}). Let us set
$K = \frac{(\alpha\beta)^{1-\delta^{\prime}}}{8C^{\prime}}$ and $L=1 + [K]$. We shall show that

\vspace{-2mm}
\begin{equation}
\label{f22}
\frac{1}{L}\sum_{1 \leq d \leq L} |{\mathcal M}(A,B,d)| \geq C^{\prime}(\alpha\beta)^{\delta^\prime} XY \; ,
\end{equation}

\noindent
so that we have $|{\mathcal M}(A,B,d)| \geq C^{\prime}(\alpha\beta)^{\delta^\prime} XY$ for some integer $d \leq L$
, which of course implies (\ref{f8}). Note that since $L$ is roughly about $\frac{(\alpha\beta)^{1-\delta^{\prime}}}{C^{\prime}}$, (\ref{f22}) is what one might expect from (\ref{f1}).

\vspace{2mm} \noindent Let $D$ be an integer in
$[\frac{2}{\alpha\beta},\frac{4}{\alpha\beta}]$. Thus
in particular $D > 1$. When $L \geq D$ we obtain (\ref{f22}) even
without (\ref{f11}). In effect, we then have $K \geq 1$ and hence
that $L < 2K$ or, what is the same thing, that $L < \frac{(\alpha
\beta)^{1-\delta^{\prime}}}{4C^{\prime}}$ from which (\ref{f22})
follows on noting that for any integer $T \geq D$, and in particular
for $T =L$, we have from (\ref{f3}) that

\vspace{-2mm}
\begin{equation}
\label{f12}  \frac{1}{T} \sum_{1 \leq d \leq T} |{\mathcal M}(A,B,d)| \; \geq \;
\left(\frac{\alpha\beta - \frac{1}{T}}{T}\right)XY \;
\geq \; \left(\frac{\alpha\beta -
\frac{1}{D}}{T}\right)XY \; \geq \; \frac{\alpha \beta XY}{2T} \; .
\end{equation}

\noindent Suppose now that $1 \leq L < D$. Let us first verify that for any integer $T$ such that $1 \leq T < D$ we have the following inequality on account of (\ref{f11}).

\vspace{-2mm}
\begin{equation}
\label{f14}
\sum_{T < d \leq D} |{\mathcal M}(A,B,d)| \; \leq \;
\left(\frac{C^{\prime}}{C}\right)^{\frac{1}{
2\delta}} (\alpha \beta)^{\frac{\delta^{\prime}}{2\delta}} T^{\frac{
1}{\delta}-1}\; (XY)^{\frac{1}{2}}
\left(\sum_{T < d \leq D} |A_d|\right)^{\frac{1}{2}}
\left(\sum_{T < d \leq D} |B_d|\right)^{\frac{1}{2}}
\;.
\end{equation}

\noindent Indeed,  for any integer $d$ satisfying $T <
d \leq D$ we have that

\vspace{-2mm}
\begin{equation}
\label{f14.5}
|A_d||B_d| = (\alpha_d X \beta_d Y)^{\frac{1}{2}} |A_d|^{\frac{1}{2}}
|B_d|^{\frac{1}{2}}\; \leq \;
\left(\frac{C^{\prime}}{C}\right)^{\frac{1}{
\delta}} (\alpha \beta)^{\frac{\delta^{\prime}}{2\delta}} T^{\frac{
1}{\delta}-1}\; (XY)^{\frac{1}{2}} |A_d|^{\frac{1}{2}}
|B_d|^{\frac{1}{2}} \; ,
\end{equation}

\noindent
where the last inequality follows from (\ref{f11}) on noting that
$d^{\frac{1}{\delta}-1} \leq T^{\frac{1}{\delta}-1}$ for $d$ satisfying $T < d \leq D$,
since $\delta \geq 1$. On combining the bound $|{\mathcal M}(A,B,d)| \leq |A_d||B_d|$ with (\ref{f14.5}) and an application of the
Cauchy-Schwarz inequality we obtain (\ref{f14}).

\vspace{2mm}
\noindent
We now estimate the sums on the right hand side of (\ref{f14}). An
application of H\"{o}lder's inequality gives

\vspace{-2mm}
\begin{equation}
\label{f15}
\sum_{T < d \leq D} |A_d| \; = \;
\sum_{T < d \leq D} \sum_{\stackrel{n \in A,}{d|n}} 1
\;\leq \; \sum_{n \in A} \tau_{D}(n) \; \leq \;
|A|^{1 -\frac{1}{q}} \left(\sum_{1 \leq n \leq X } \tau_{D}(n)^{q}\right)^{\frac{1}{q}} \;.
\end{equation}

\noindent From Lemma 2.1 we have the upper bound $c(q)DX$ for the
last sum in (\ref{f15}). Since $|A| \leq  2\alpha X$
and $D \leq \frac{4}{\alpha \beta}$, we deduce from (\ref{f15}) that

\vspace{-2mm}
\begin{equation}
\label{f15.1} \sum_{T < d \leq D} |A_d| \leq
(2\alpha)^{1-\frac{2}{q}}\beta^{-\frac{1}{q}}(4c(q))^{\frac{1}{q}}X.
\end{equation}

\noindent
Arguing similarly, we obtain the bound

\vspace{-2mm}
\begin{equation}
\label{f15.2} \sum_{T < d \leq D} |B_d| \leq
(2\beta)^{1-\frac{2}{q}}\alpha^{-\frac{1}{q}}(4c(q))^{\frac{1}{q}}Y.
\end{equation}

\noindent
With these estimates we conclude from (\ref{f14}) that for any integer $T$ satisfying $1 \leq T < D$ we have

\vspace{-2mm}
\begin{equation}
\label{f16} \sum_{T < d \leq D} |{\mathcal M}(A,B,d)| \; \leq \;
2\left(\frac{C^{\prime}}{C}\right)^{ \frac{1}{2\delta}}
(\alpha \beta)^{ \frac{\delta^{\prime}}{2\delta} + \frac{1}{2}
-\frac{3}{2q}} T^{\frac{1}{\delta} -1} (4c(q))^{\frac{1}{q}} XY \; ,
\end{equation}

\noindent We now reveal that our choices for $C^{\prime}$ and
$\delta^{\prime}$ were  made so that $K$ satisfies the relation

\vspace{-2mm}
\begin{equation}
\label{f17} 2\left(\frac{C^{\prime}}{C}\right)^{\frac{1}{2\delta}}
(\alpha \beta)^{ \frac{\delta^{\prime}}{2\delta} + \frac{1}{2}
-\frac{3}{2q}} K^{\frac{1}{\delta} -1} (4c(q))^{\frac{1}{q}} \; = \;
\frac{\alpha\beta}{4} \; ,
\end{equation}

\noindent as may be confirmed by a modest calculation
using the expressions defining $C^{\prime}$ and $\delta^{\prime}$ in
terms of $C$ and $\delta$.

\vspace{2mm} \noindent We see that $\sum_{L < d \leq D} |{\mathcal
M}(A,B,d)| \leq \frac{\alpha \beta}{ 4}
XY$ using (\ref{f16}) for $T = L$ together with (\ref{f17}) and
noting that $K < L$. Since (\ref{f12}) applied with $T =D$ gives us
$\sum_{1 \leq d \leq D} |{\mathcal M}(A,B,d)| \geq
\frac{\alpha \beta}{2} XY$, we conclude that when $1
\leq L < D$ we have

\vspace{-2mm}
\begin{equation}
\label{f18}
\frac{1}{L}\sum_{1 \leq d \leq L} |{\mathcal M}(A,B,d)| \, \geq \, \frac{\alpha \beta}{4L} XY \; .
\end{equation}

\vspace{2mm}
\noindent
If $L =1$ we obtain (\ref{f22}) from (\ref{f18}) on noting that $\frac{\alpha \beta}{4} \geq C^{\prime} (\alpha \beta)^{\delta^{\prime}}$, since $\frac{1}{4} \geq C^{\prime}$ and $1 \leq \delta^{\prime} $. When $1 < L < D$ we have $1 \leq K$ and hence $L < \frac{(\alpha \beta)^{1-\delta^{\prime}}}{4C^{\prime}}$ so that (\ref{f22}) results from (\ref{f18}) in this final case as well.

\vspace{3mm}
\begin{cor} Every $\delta>1$ is an admissible exponent. \end{cor}

\vspace{1mm}
\noindent
{\sc Proof. ---} Let $q$ be any integer $\geq 4$ and let $\{\delta_n (q)\}_{n \geq 1}$ the sequence
of real numbers determined by the relations $\delta_1 (q) = 2$ and

\vspace{-2mm}
\begin{equation}
\label{fe0}
\delta_{n+1} (q) = \frac{3\delta_n(q)\left (1+\frac 1q\right
)-2}{2\delta_n(q)-1}
\end{equation}

\noindent
 for $n \geq 1$. Then each $\delta_n(q)$ is
 an admissible exponent by Propositions 2.1 and 2.2. It is easily verified that the sequence $\delta_n(q)$ is decreasing and has a limit $\delta(q)$ given by the relation

\vspace{-2mm}
\begin{equation}
\label{fe1}
\delta(q)=1+\frac{3}{4q}+\frac 12\sqrt{\frac
 6q+\frac{9}{4q^2}}.
\end{equation}

\noindent
Plainly, any $\delta>\delta(q)$ is an admissible exponent. The corollary now follows on taking $q$ arbitrarily large in (\ref{fe1}).

\vspace{2mm}
\noindent
Theorem 1.1 follows from the above corollary and the definition of admissible exponents on recalling that $|A/B| \geq \sup_{d \geq 1} |{\mathcal M}(A,B,d)|$.

\end{nsc}

\begin{nsc}{Counterexamples}

\noindent
Let us first prove Theorem 1.2. To this end, given an integer $m \geq 1$ let ${\mathcal P}$ denote
any set of $2m$ prime numbers and, for any subset $I$ of ${\mathcal P}$, let $d(I) = \prod_{p \in I} p$.
If $S({\mathcal P})$ denotes the set of all $d(I)$ with $|I| = m$, we have the following lemma.

\vspace{2mm}
\begin{le}
For any $\epsilon >0$, we have $|S({\mathcal P})/S({\mathcal P})| \leq \epsilon |S({\mathcal P})|^2$ for all sufficiently large $m$.
\end{le}

\vspace{2mm}
\noindent
{\sc Proof. ---} Plainly, we have $|S({\mathcal P})| = \binom{2m}{m}$. Let ${\mathcal Q}$ be the set of
ordered pairs of disjoint subsets of ${\mathcal P}$. Then, for any $I$ and $J$ subsets of ${\mathcal P}$, we have

\vspace{-2mm}
\begin{equation}
\label{e1}
\frac{d(I)}{d(J)} = \frac{d(I\setminus J)}{d(J\setminus I)} \; ,
\end{equation}

\noindent
and since $I\setminus J$ and $J\setminus I$ are disjoint, $(I\setminus J,J \setminus I)$ is in ${\mathcal Q}$.
Thus $|S({\mathcal P})/S({\mathcal P})| \leq |{\mathcal Q}|$. Let us associate any $(U,V)$ in ${\mathcal Q}$ to the map from ${\mathcal P}$ to the three element set $\{1,2,3\}$ that takes $U$ to 1, $V$ to 2 and the complement of $U \cup V$ in ${\mathcal P}$ to 3. It is easily seen that this association in fact gives a bijection from ${\mathcal Q}$ onto the set of maps from ${\mathcal P}$ to $\{1,2,3\}$ and hence that $|{\mathcal Q}| = 3^{2m}$. In summary, we deduce that

\vspace{-2mm}
\begin{equation}
\label{s31}
|S({\mathcal P})/S({\mathcal P})| \; \leq \; |{\mathcal Q}|\; =\; 3^{2m} \; = \; \frac{3^{2m}}{{\binom{2m}{m}}^2}\,|S({\mathcal P})|^2 \; \leq \; (2m+1)^2 \left(\frac{3}{4}\right)^{2m}|S({\mathcal P})|^2 \; ,
\end{equation}

\noindent
where we have used the inequality $\binom{2m}{m} \geq \frac{2^{2m}}{2m+1}$.
The lemma follows from (\ref{s31}) on noting that $(2m+1)^2 \left(\frac{3}{4}\right)^{2m} \rightarrow 0$ as $m \rightarrow +\infty$.

\vspace{2mm} \noindent
{\sc Proof of Theorem 1.2. ---} Given an integer $m \geq
1$, it is easily deduced from the prime number theorem that the
interval $[T, T + T/m]$ contains at least $2m$ prime numbers when
$T$ is sufficiently large. For such a $T$, let ${\mathcal P}$ be a
subset of $2m$ prime numbers in $[T, T + T/m]$. If ${\mathcal
A}({\mathcal P})$ is the sequence of integers $\geq 1$ that are
divisible by at least one of the integers $d(I)$ in $S({\mathcal
P})$ then a simple application of the principle of inclusion and
exclusion implies that ${\mathcal A}({\mathcal P})$ has an
asymptotic density $\alpha({\mathcal P})$ that is given by the
relation

\vspace{-2mm}
\begin{equation}
\label{asym}
\alpha({\mathcal P}) \; = \;
\sum_{1 \leq r \leq \binom{2m}{m}} (-1)^{r+1} \sum_{1\leq i_1 < i_2 \ldots < i_r \leq \binom{2m}{m}} \frac{1}{d(I_{i_1} \cup I_{i_2} \ldots \cup I_{i_r})} \; ,
\end{equation}

\vspace{2mm}
\noindent
where $I_1,I_2, \ldots, I_{\binom{2m}{m}}$ are the subsets of cardinality $m$ in ${\mathcal P}$.

\vspace{2mm}
\noindent
For any $i$ we have $T^m \leq
d(I_i) \leq (1+\frac 1m)^m T^{m}<eT^m$. Consequently, for the term
$r =1$ in (\ref{asym}) we obtain

\vspace{-2mm}
\begin{equation}
\sum_{1\leq i \leq \binom{2m}{m}} \frac{1}{d(I_i)}\; \geq \;
\frac{\binom{2m}{m}}{eT^m}.
\end{equation}

\noindent When $r \geq 2$, we have that $d(I_{i_1} \cup I_{i_2}
\ldots \cup I_{i_r})$, for any distinct indices $i_1, i_2 \ldots,
i_r$, has at least $k+1$ prime factors in ${\mathcal
P}$ and hence is $\geq T^{m+1}$. It follows from (\ref{asym}) and
these bounds that we have

\vspace{-2mm}
\begin{equation}
\label{e2.1} \alpha({\mathcal P})\ge
\frac{\binom{2m}{m}}{eT^m}-\frac{2^{\binom{2m}{m} }}{T^{m+1}}\ge
\frac{\binom{2m}{m}}{3T^m}
\end{equation}

\noindent when $T$ is sufficiently large. In particular, on
recalling that $|S({\mathcal P})| = \binom{2m}{m}$, we obtain that
for any integer $m \geq 1$, we have

\vspace{-2mm}
\begin{equation}
\label{e2} \alpha({\mathcal P}) \; \geq \; \frac{|S({\mathcal
P})|}{3T^m} \;
\end{equation}

\noindent for all sufficiently large $T$ and ${\mathcal P}$ any set
of $2m$ prime numbers in $[T, T+T/m]$.

\vspace{2mm}
\noindent
Finally, for ${\mathcal P}$ as above and any $X \geq 1$, let us set
 $A = {\mathcal A}({\mathcal P}) \cap [1,X]$.
 Since $\alpha({\mathcal P})$ is the asymptotic density of ${\mathcal A}({\mathcal P})$, we
 have from (\ref{e2}) that $|A| \geq \frac{|S({\mathcal P})|}{4T^m}X$, for all large enough $X$ and $T$. Clearly, each integer in $A$ is of the form $d(I)n$, for some $d(I)$ in $S({\mathcal P})$ and an integer $n$, which must necessarily be $\leq \frac{X}{T^m}$, since $A$ is in $[1,X]$ and $d(I) \geq T^m$. Consequently, we have we have $|A/A| \leq \frac{|S({\mathcal P})/S({\mathcal P})|}{T^{2m}} X^2$, for all large enough $X$ and $T$. On comparing $|A|$ and $|A/A|$ by means of Lemma 3.1, we see that $A$ meets the conditions of Theorem 1.2 when $m$, $T$ and $X$ are all sufficiently large.

\vspace{3mm} \noindent  {\sc Proof of Theorem 1.3. ---}
The number of primitive integer points, that is, integer points with
coprime co-ordinates, in $[1,\gamma X] \times [1,\gamma Y]$ is $\sim
\frac{6}{\pi^2} \gamma^2 XY$ as $X$, $Y \rightarrow \infty$. Thus
for any $\gamma$ in $(0,1]$ and all sufficiently large $X$ and $Y$,
there is a subset $S$ of the primitive integer points in $[1,\gamma
X] \times [1,\gamma Y]$ satisfying $\frac{\gamma^2}{4} XY \leq |S|
\leq \frac{\gamma^2}{2} XY$. Let us take for $C$ the union of the
sets $d.S$ with $d$ varying over the interval
$[1,\frac{1}{\gamma}]$, where each $d.S$ is the set of $(da,db)$
with $(a,b)$ varying over $S$. Then $C$ is contained in $[1, X] \times [1,
Y]$. Moreover, the sets $d.S$ are disjoint but ${\rm Frac}(d.S) =
{\rm Frac}(S)$, for each $d$, and 
$|{\rm Frac}(S)| = |S|$. We therefore have $|C| = [\frac{1}{\gamma}]|S|
\geq \frac{\gamma}{8} XY$ but $|{\rm Frac}(C)|= |{\rm Frac}(S)| =
|S|\leq \frac{\gamma^2}{2} XY$.

\end{nsc}

\begin{nsc}{Gaps in Product Sequences}

\noindent
We now deduce Theorem 1.4 from Theorem 1.1. Let ${\mathcal A}$ and ${\mathcal B}$ be sequences with upper asymptotic densities $\alpha$ and $\beta$.
Then there exist infinitely many real numbers $X$ and $Y$ $\geq 1$ such that $|{\mathcal A} \cap (\frac{X}{2},X]| \geq \frac{\alpha X}{4}$ and
 $|{\mathcal B} \cap (\frac{Y}{2},Y]| \geq \frac{\beta Y}{4}$. For such $X$ and $Y$ let us apply Theorem 1.1  to the sets $A = {\mathcal A} \cap (\frac{X}{2},X]$
 and $B = {\mathcal B} \cap (\frac{Y}{2},Y]$. We then have that $|A/B| \gg (\alpha \beta)^{1+\epsilon} XY$, where the implied constant
 depends on $\epsilon$ alone. Since $A/B$ is a subset of the interval $[\frac{X}{2Y}, \frac{2X}{Y}]$, which is of length $\frac{X}{Y}$,
 we deduce that there are distinct $a/b$ and $a^{\prime}/b^{\prime}$ in $A/B$ such that

\vspace{-2mm}
\begin{equation}
\label{e5} 0< \left|\frac{a}{b} -
\frac{a^{\prime}}{b^{\prime}}\right| \ll \frac{X/Y}{(\alpha
\beta)^{1+\epsilon}XY}  = \frac{1}{(\alpha \beta)^{1+\epsilon}Y^2}
\; .
\end{equation}

\noindent
Since $|bb^{\prime}| \leq Y^2$, it
 follows from (\ref{e5}) that difference between the
 distinct terms $ba^{\prime}$ and $b^{\prime}a$ of the product
 sequence ${\mathcal A}.{\mathcal B}$ is $\ll \frac{1}{(\alpha \beta)^{1+\epsilon}}$.
 Since there are infinitely many distinct $X$ and $Y$ satisfying the required conditions,
 there are infinitely many such pairs of terms in the product sequence ${\mathcal A}.{\mathcal B}$.

\end{nsc}

\vspace{3mm}
\noindent
{\bf Acknowledgments : }
We arrived at Theorem 1.3 in response to a question of Professor Adrian Ubis, whom we gladly thank. We also wish to thank
the CRM, Barcelona and HRI, Allahabad for
opportunities that supported discussions on the problems addressed here.

\vspace{2mm} \noindent J. Cilleruelo was supported by
Grant CCG08-UAM/ESP-3906 and MTM2008-03880 of the MYCIT
, Spain during the course of preparation of this note. D.S. Ramana
is with the Harish-Chandra Research Institute which is a
constituent institution of the Homi Bhabha National Institute,
Mumbai, India. Olivier Ramar\'{e} is supported by the CNRS, France.


\vspace{-2mm}
\begin{thebibliography}{1}

\bibitem{ber}
G.~Berczi.
\newblock On the distribution of products of members of a sequence
with positive density,
\newblock {\em Per. Math. Hung.}, 44(2002):137-145.

\bibitem{ch}
J.~Cilleruelo and T.H.~Le,
\newblock Gaps in product sequences.
\newblock {\em Israel Journal of Mathematics}, {to
appear.}

\bibitem{san}
C. Sandor,
\newblock On the minimal gaps between products of members of a sequence with positive density.
\newblock {\em Ann. Univ. Sci. Budapest E\"{o}tv\"{o}s Sect. Math.} 28(2005), 3-7.

\bibitem{sar}
A. Sark\"{o}zy,
\newblock Unsolved problems in number theory,
\newblock {\em Per. Math. Hung.}, 42(2001), 17-36.

\end{thebibliography}
\end{document}